\documentclass{article}

\usepackage{arxiv}

\usepackage[utf8]{inputenc} 
\usepackage[T1]{fontenc}    
\usepackage{hyperref}       
\usepackage{url}            
\usepackage{booktabs}       
\usepackage{amsfonts}       
\usepackage{nicefrac}       
\usepackage{microtype}      
\usepackage{lipsum}

\usepackage{tikz}

\usetikzlibrary{matrix,arrows,shapes.geometric}

\usepackage[all,2cell,cmtip]{xy}

\usepackage{mathtools}

\usepackage{amssymb}

\usepackage{enumerate}

\usepackage[all]{xy}
\xyoption{2cell}
\xyoption{curve}
\UseTwocells
\input{diagxy}
\CompileMatrices      
\usepackage{tikz}
\usepackage{pdfpages}

\usepackage{mathrsfs}

\usepackage{tikz}

\usepackage{stmaryrd}

\usepackage{scalerel,amssymb} 
\def\msquare{\mathord{\scalerel*{\Box}{gX}}} 

\usepackage{amsthm}

\newtheorem{lema}{Lemma}

\newtheorem{dfn}{Definition}

\newtheorem{thm}{Theorem}

\newtheorem{pps}{Proposition}

\newtheorem{cor}{Corollary}

\newtheorem{rem}{Remark}

\title{Implications of a Quillen Model Structures-Based Framework for Locality under Logical Equivalence}

\author{
  Hendrick Maia
  \\
  \texttt{hendrickmaia@gmail.com} \\
}

\begin{document}
\maketitle

\begin{abstract}

In [15] a homotopic variation for locality of logics was presented, namely a Quillen model category-based framework for locality under logical equivalence, for every primitive-positive sentence of quantifier-rank $k$. In this paper, we will present some of the implications and possible themes for investigations that arise from the aforementioned framework. 

\end{abstract}


\section{Introduction}

Locality is a property of logics, whose origins lie in the works of Hanf [10] and Gaifman [8], having their utility in the context of finite model theory. Such a property is quite useful in proofs of inexpressibility, but it is also useful in establishing normal forms for logical formulas. 

There are generally two forms of locality: (i') if two structures $\mathfrak{A}$ and $\mathfrak{B}$ realize the same multiset of types of neighborhoods of radius $d$, then they agree on a given sentence $\Phi$. Here $d$ depends only on $\Phi$; (ii')  if the $d$-neighborhoods of two tuples $\vec{a}_1$ and $\vec{a}_2$ in a structure $\mathfrak{A}$ are isomorphic, then $\mathfrak{A} \models \Phi(\vec{a}_1) \Leftrightarrow \Phi(\vec{a}_2)$. Again, $d$ depends on $\Phi$, and not on $\mathfrak{A}$. Form (i') originated from Hanf's works [10]. Form (ii') came from Gaifman's theorem [8]. Before proceeding, I will establish some notation. 

\noindent \textbf{Notations:} All structures here are finite, whose vocabularies are finite sequences of relation symbols $\sigma = \langle R_1,...,R_l \rangle$. A $\sigma$-structure $\mathfrak{A}$ consists of a finite universe $A$ and an interpretation of each $p_i$-ary relation symbol $R_i$ in $\sigma$ as $R^{\mathfrak{A}}_i \subseteq A^{p_i}$. Isomorphism of structures will be denoted by $\cong$. I shall use the notation $\sigma_n$ for $\sigma$ expanded with $n$ constant symbols. 

The quantifier-rank of a formula $\Phi$ is the maximal nesting depth of quantifiers in $\Phi$. 

Given a structure $\mathfrak{A}$, its \emph{Gaifman graph} $\mathcal{G}(\mathfrak{A})$ is defined as $\langle A,E \rangle$ where $(a,b)$ is in $E$ if, and only if there is a tuple $\vec{c} \in R^{\mathfrak{A}}_i$ for some $i$ such that both $a$ and $b$ are in $\vec{c}$. The distance $d(a,b)$ is defined as the length of the shortest path from $a$ to $b$ in $\mathcal{G}(\mathfrak{A})$; we assume $d(a,a) = 0$. If $\vec{a} = (a_1,...,a_n)$, then $d(\vec{a},b) = \mathrm{min}_id(a_i,b)$. Given $\vec{a}$ over $A$, its \emph{$r$-ball} $B_r^{\mathfrak{A}}(\vec{a})$ is $\{b \in A \mid d(\vec{a},b) \leq r\}$. If $|\vec{a}| = n$, its \emph{$r$-neighborhood} $N_r^{\mathfrak{A}}(\vec{a})$ is defined as a $\sigma_n$-structure
$$\langle B_r^{\mathfrak{A}}(\vec{a}), R_1^{\mathfrak{A}} \cap B_r^{\mathfrak{A}}(\vec{a})^{p_1},..., R_l^{\mathfrak{A}} \cap B_r^{\mathfrak{A}}(\vec{a})^{p_l}, a_1,...,a_n \rangle.$$
Note that for any isomorphism $h:N_r^{\mathfrak{A}}(\vec{a}) \rightarrow N_r^{\mathfrak{B}}(\vec{b})$ it must be the case that $h(\vec{a}) = \vec{b}$.

Given a tuple  $\vec{a} = (a_1,...,a_n)$ and an element $c$, we write  $\vec{a}c$ for the tuple $(a_1,...,a_n,c)$. 

An $m$-ary query, $m \geq 0$, on $\sigma$-structures, is a mapping $Q$ that associates with each structure $\mathfrak{A}$ a subset of $A^m$, such that $Q$ is closed under isomorphism: if $\mathfrak{A} \cong \mathfrak{B}$ via isomorphism $h : A \rightarrow B$, then $Q(\mathfrak{B}) = h(Q(\mathfrak{A}))$.

We write $\mathfrak{A} \equiv_k \mathfrak{B}$ if $\mathfrak{A}$ and $\mathfrak{B}$ agree on all FO sentences of quantifier rank up to $k$, and $(\mathfrak{A},\vec{a}) \equiv_k (\mathfrak{B},\vec{b})$ if $\mathfrak{A} \models \Phi(\vec{a}) \Leftrightarrow \mathfrak{B} \models \Phi(\vec{b})$ for every FO formula $\Phi(\vec{x})$ of quantifier rank up to $k$ ($k$-logical equivalence). It is well known (see [14]) that $\mathfrak{A} \equiv_k \mathfrak{B}$ if, and only if, the duplicator has a winning strategy in the $k$-round Ehrenfeucht-Fra\"iss\'e game on $\mathfrak{A}$ and $\mathfrak{B}$, and $(\mathfrak{A},\vec{a}) \equiv_k (\mathfrak{B},\vec{b})$ if, and only if, the duplicator has a winning strategy in the $k$-round Ehrenfeucht-Fra\"iss\'e game on $\mathfrak{A}$ and $\mathfrak{B}$ starting in position $(\vec{a},\vec{b})$.

The tree-depth $\mathrm{td}(\mathscr{G})$ of a finite graph $\mathscr{G}$ is the minimal height of a finite rooted forest whose closure contains $\mathscr{G}$ as a subgraph. The tree-depth $\mathrm{td}(\mathfrak{A})$ of a finite structure $\mathfrak{A}$ is defined as the tree-depth of the Gaifman graph of $\mathfrak{A}$.

Let $k \in \mathbb{N}$. We say $\mathfrak{A}$ is \emph{$k$-homomorphic} to $\mathfrak{B}$ if $\mathfrak{C}$ is homomorphic to $\mathfrak{A}$ and $\mathfrak{C}$ is homomorphic to $\mathfrak{B}$ for every finite structure $\mathfrak{C}$ of tree-depth at most $k$. We say $\mathfrak{A}$ and $\mathfrak{B}$ are \emph{$k$-homomorphically equivalent} if $\mathfrak{A}$ is $k$-homomorphic to $\mathfrak{B}$ and $\mathfrak{B}$ is $k$-homomorphic to $\mathfrak{A}$.  

A $\sigma$-structure $\mathfrak{C}$ is called a core if $\mathfrak{C}$ is not homomorphic to any proper substructure of $\mathfrak{C}$. Let $k \in \mathbb{N}$. A $k$-core is a core with tree-depth at most $k$.

\qed

There is no doubt about the usefulness of the notion of locality, which as seen applies to a huge number of situations. However, there is a deficiency in such a notion: all versions of the notion of locality refer to isomorphism of neighborhoods, which is a fairly strong property. For example, where structures simply do not have sufficient isomorphic neighborhoods, versions of the notion of locality obviously cannot be applied. So the question that immediately arises is: would it be possible to weaken such a condition, and maintain Hanf / Gaifman-localities?

Arenas, Barceló and Libkin [1] establish a new condition for the notions of locality, weakening the requirement that neighborhoods should be isomorphic, establishing only the condition that they must be indistinguishable in a given logic. That is, instead of requiring $N_d(\vec{a}) \cong N_d(\vec{b})$, you should only require $N_d(\vec{a}) \equiv_k N_d(\vec{b})$, for some $k \geq 0$. Using the fact that logical equivalence is often captured by Ehrenfeucht – Fra{\"i}ssé games, the authors formulate a game-based framework in which logical equivalence-based locality can be defined. Thus, the notion defined by the authors is that of \emph{game-based locality}.

Note that the intuitive point from which the authors start is the idea of \emph{neighborhood indistinguishability}. Thus, the intuition behind the notion of game-based locality is to describe the indistinguishability of neighborhoods in terms of \emph{winning game strategies}. To achieve the necessary generalization, Arenas, Barceló and Libkin define an abstract view of the games that characterize the expressiveness of logics that are local under isomorphism. The basic idea is as follows: in each round the duplicator has a set of functions (tactics) that will determine his responses to possible moves by the spoiler. 

The problem with the game-based framework, which can be seen as a general problem, can be found by thinking about why the notion of locality has gained so much space. This is because winning games are nontrivial, even for very simple examples. That is, even for fairly easy examples, the difficulty of winning games is quite high. Thus, for its simplicity, the notion of locality eventually gained much attention, as well as further developments and extensions beyond first-order logic (see [14]). However, the need to weaken the notion of locality brought back precisely what the notion of locality avoided, namely games. So why should we go back to working with complicated game methods? We are using the notion of locality exactly to avoid games! Therefore, a game-based framework for weakening the notion of locality does not seem to be very plausible. In the following, I will outline the three specific problems that the game-based framework has.

The question that immediately arises is: is it possible to define the notion of locality under logical equivalence without resorting to game-based frameworks? 

This question motivated the development of a model category-based framework for locality under logical equivalence. Thus, in [15] the following result was stated and proved:

\begin{thm}[MAIA]
	There is a Quillen model structure $\mathbb{M}$ on $\mathbf{STRUCT}[\sigma_n]_{(d,0)}^{\widetilde{\mathcal{T}(\sigma_n)}}$ such that the homotopic equivalences in $\mathbb{M}$ coincides with the homomorphic equivalences in $\mathbf{STRUCT}[\sigma_n]_{(d,0)}^{\widetilde{\mathcal{T}(\sigma_n)}}$, and such that for every $k$-homotopic equivalence, $\sim_k$, and every $k$-logical equivalence, $\equiv_k$, $N_{d}^{\mathfrak{A}}(\vec{a}) \sim_k N_{d}^{\mathfrak{B}}(\vec{b})$ if and only if $N_{d}^{\mathfrak{A}} (\vec{a}) \equiv_k N_{d}^{\mathfrak{B}}(\vec{b})$, for every primitive-positive sentence \footnote{This is merely a matter of convenience; Theorem 1 remains valid when stated more generally for formulas instead of sentences.} with quantifier-rank $k$.
\end{thm}

The category $\mathbf{STRUCT}[\sigma_n]_{(d,0)}^{\widetilde{\mathcal{T}(\sigma_n)}}$ (see [15] \S 4) is the category $\mathbf{STRUCT}[\sigma_n]_d$ of $d$-neighborhoods and homomorphisms between such $d$-neighborhoods that admits 0-neighborhoods as objects, and which has the $\sigma_n$-structure $\widetilde{\mathcal {T}(\sigma_n)}$ (see [15] Definition 23-24) as an object.

 Every object is fibrant and cofibrant in the model structure $\mathbb{M}$. The acyclic fibrations are exactly the retractions in $\mathbf{STRUCT}[\sigma_n]_{(d,0)}^{\widetilde{\mathcal{T}(\sigma_n)}}$. Any two morphisms with same domain and codomain are homotopic in the model structure $\mathbb{M}$. Weak equivalences in $\mathbb{M}$ are the maps inducing isomorphisms on the cores, i.e., homomorphic equivalences in $\mathbf{STRUCT}[\sigma_n]_{(d,0)}^{{\widetilde{\mathcal{T}(\sigma_n)}}}$, which in $\mathbb{M}$ coincide with homotopic equivalences. The homotopy category of the model structure $\mathbb{M}$ is $\mathbf{STRUCT}[\sigma_n]_{(d,0)}^{{\widetilde{\mathcal{T}(\sigma_n)}}^{\mathrm{core}}}$.

Theorem 1 naturally extends when we consider the category of finite $\sigma$-structures and homomorphisms, $\mathbf{STRUCT}[\sigma]$, since it is easy to see that $\mathbf{STRUCT}[\sigma]$ is finitely complete and cocomplete. So, we also have

\begin{thm}
	There is a Quillen model structure $\mathbb{M}$ on $\mathbf{STRUCT}[\sigma]$ such that the homotopic equivalences in $\mathbb{M}$ coincides with the homomorphic equivalences in $\mathbf{STRUCT}[\sigma]$, and such that for every $k$-homotopic equivalence, $\sim_k$, and every $k$-logical equivalence, $\equiv_k$, $\mathfrak{A} \sim_k \mathfrak{B}$ if and only if $\mathfrak{A} \equiv_k \mathfrak{B}$, for every primitive-positive sentence with quantifier-rank $k$.
\end{thm}

Although Theorem 1 (and Theorem 2) remains valid only for primitive-positive sentences, it is valid for all sentences if we consider only a special class of structures.

The notation $(\mathfrak{A},\vec{a}) \rightarrow_X (\mathfrak{B},\vec{b})$ (introduced in [19] \S 2.2) to express that there exists a homomorphism from $\mathfrak{A}$ to $\mathfrak{B}$ over $X$ which carries tuple $\vec{a}$ to tuple $\vec{b}$. The notation extends to $k$-homomorphism over $X$ in the obvious way.

\begin{dfn}
A structure $\mathfrak{A}$ is $k$-extendable if, for every set $X \subseteq A$ of size $< k$ and every structure $\mathfrak{B}$ such that $\mathfrak{A} \rightleftarrows_X^{k-|X|} \mathfrak{B}$, it holds that $\forall b \in B \textrm{ } \exists a \in A$ s.t. $\mathfrak{A} \rightleftarrows_X^{k-|X|-1} \mathfrak{B}$.
\end{dfn}{}

\begin{lema}
Suppose structure $\mathfrak{A}$ and $\mathfrak{B}$ are $k$-extendable and $\mathfrak{A} \rightleftarrows^n \mathfrak{B}$. Then $\mathfrak{A} \equiv_k \mathfrak{B}$.
\end{lema}{}

\begin{proof}
([19], p. 30).
\end{proof}{}

\begin{cor}
For $k$-extendables $\sigma_n$-structures, Theorem 1 (and Theorem 2) holds for every sentence with quantifier-rank $k$.
\end{cor}{}

In this paper, we will present some of the implications and possible themes for investigations that arise from Theorem 1 and Theorem 2.

{\bf Overview of the paper:} In {\bf Section 2} we present the implications of Theorem 1 for the theme of locality. In short, Section 2 can be seen as a homotopic variation for the locality of logics. Of particular interest are Definitions 8 and 9, which takes us to the concept of localization, as it usually comes up in topology.  {\bf Section 3}  we briefly present how the homotopy category of any model category is naturally a $\mathrm{Ho(SSet)}$-module category. By Theorems 1 and 2, this obviously implies that the homotopy categories of $\mathbf{STRUCT}[\sigma_n]_{(d,0)}^{\widetilde{\mathcal{T}(\sigma_n)}}$ and $\mathbf{STRUCT}[\sigma]$ are naturally $\mathrm{Ho(SSet)}$-module categories. The main point of this section is the fact that the homotopy category of a pointed model category is a pre-triangulated category. The paper ends in {\bf Section 4}, where we briefly show that there are several specific results when considering $\mathrm{Ho(SSet_*)}$-module categories, that is, when we are working within the scope of homotopy categories of pointed model categories. Then, we show that with under and over categories, it is possible to obtain pointed model categories from categories $\mathbf{STRUCT}[\sigma_n]_{(d,0)}^{\widetilde{\mathcal{T}(\sigma_n)}}$ and $\mathbf{STRUCT}[\sigma]$. This shows that it is possible to investigate how we can apply the aforementioned specific results within the scope of finite structures. In particular, we have that the homotopy categories of the pointed model categories built from $\mathbf{STRUCT}[\sigma_n]_{(d,0)}^{\widetilde{\mathcal{T}(\sigma_n)}}$ and $\mathbf{STRUCT}[\sigma]$ are pre-triangulated categories.

\section{A Homotopic Variation on the Theme of Locality}
\label{sec:others}

Theorem 1 allows you to define locality under logical equivalence without game-based frameworks. That is, different from what happens with the approach of Arenas, Barceló and Libkin, who start from a game-based framework (game-based locality), that is, describe the logical indistinguishability of neighborhoods in terms of $\mathfrak{F}$-games, the approach proposed here is that of a Quillen model categories-based framework (locality under $k$-homotopic equivalence, for some $k$ [15]), that is, the purpose here is to describe logical indistinguishability of neighborhoods in terms of homotopic notions. This is interesting not only because it is an alternative to the game-based framework, but also because it opens up a new range of possibilities for working with locality under logical equivalence, namely the whole technical apparatus that comes up with Quillen model categories.

Let $\mathrm{STRUCT}[\sigma_n]$ be the class of $\sigma_n$-structures. First, note that it is possible to define a $k$-homotopic version of $d$-equivalence, which I define as follows:

\begin{dfn}
	Consider $\vec{a} \in A^n $ and $\vec{b} \in B^n$. So two $\sigma_n$-structures $\mathfrak{A}, \mathfrak{B}$ are $d$-equivalents under $k$-homotopic equivalence if, and only if, there is a bijection $f : A \rightarrow B$ such that $N_d^{\mathfrak{A}}(\vec{a} c) \sim_k N_d^{\mathfrak{B}}(\vec{b} f(c))$, for every $c \in A$. When this occurs, I denote this fact by $(\mathfrak{A},\vec{a}) \leftrightarrows_{(d,\sim_k)} (\mathfrak{B}, \vec{b})$. If $n = 0$, we have to satisfy only the condition that there is a bijection $f : A \rightarrow B$ such that $N_d^{\mathfrak{A}}(c) \sim_k N_d^{\mathfrak{B}}(f(c))$, for every $c \in A$, and the denotation is simply $\mathfrak{A} \leftrightarrows_{(d, \sim_k)} \mathfrak{B}$.
\end{dfn}

The notions of Hanf/Gaifman-localities can be defined under $k$-homotopic equivalence rather than under isomorphisms:

\begin{dfn}[Gaifman-locality under $\sim_k$-equivalence]
	
	An $m$-ary query $Q$, $m > 0$, on $\sigma_n$-structures, is called \emph{Gaifman-local under $\sim_k$-equivalence} if there exists a number $d \geq 0$ such that for every $\sigma$-structure $\mathfrak{A}$ and every $\vec{a}_1, \vec{a}_2 \in A^m$,
	
	$$N_d^{\mathfrak{A}}(\vec{a}_1) \sim_k N_d^{\mathfrak{A}}(\vec{a}_2) \textrm{ implies } (\vec{a}_1 \in Q(\mathfrak{A}) \Leftrightarrow \vec{a}_2 \in Q(\mathfrak{A})).$$
	
	\noindent The minimum $d$ for which the above condition holds is called the \emph{Gaifman-locality rank of $Q$ under $\sim_k$-equivalence}, and is denoted by $\mathsf{lr}_{\sim_k}(Q)$.

\end{dfn}

\begin{dfn}[Hanf-locality under $\sim_k$-equivalence]
	
	An $m$-ary query $Q$ on $\sigma$-structures is  \emph{Hanf-local under $\sim_k$-equivalence} if there exists a number $d \geq 0$ such that for every $\mathfrak{A}, \mathfrak{B} \in \mathrm{STRUCT}[\sigma]$,
	
	$$(\mathfrak{A},\vec{a}) \leftrightarrows_{(d,\sim_k)} (\mathfrak{B},\vec{b}) \textrm{ implies } (\vec{a} \in Q(\mathfrak{A}) \Leftrightarrow \vec{b} \in Q(\mathfrak{B})).$$
	
	\noindent The smallest $d$ for which the above condition holds is called the \emph{Hanf-locality rank of $Q$ under $\sim_k$-equivalence}, and is denoted by $\mathsf{hlr}_{\sim_k}(Q)$.
	
\end{dfn}

By Theorem 1, Hanf/Gaifman-localities under $k$-homotopic equivalence coincide with Hanf/Gaifman-localities under $k$-logical equivalence, for every primitive-positive sentence of quantifier-rank $k$. Now, I am going to show a more interesting implication of Theorem 1.

There are contexts in which we have a category $\mathcal{C}$ that is misbehaving (in a context-dependent sense, of course): it may be that $\mathcal{C}$ does not have any desired properties. One solution to such a problem is to try to find a second category, $\hat{\mathcal{C}}$, that has the same objects as $\mathcal{C}$, along with a functor $\mathcal{C} \rightarrow \hat{\mathcal{C}}$, which is the identity about objects, and such that $\hat{\mathcal{C}}$ (in a context-dependent sense, of course) is better behaved than $\mathcal{C}$, while it can be considered as an approximation of $\mathcal{C}$ (also in a context-dependent sense, of course). Sure, some $\mathcal{C}$ structure may be lost along the way, but that may be a small price to pay when $\hat{\mathcal{C}}$ shows us new insights and solutions that $\mathcal{C}$ is not able to provide.

What if we had a non-game-based framework for logic-based locality that allowed us, in some sense, to recover the isomorphic indistinguishability of neighborhoods? In other words, what if we had a ''approximation'' of $\mathbf{STRUCT}[\sigma_n]_{(d, 0)}^{\widetilde{\mathcal{T}(\sigma_n)}}$ that would allow us to treat logical indistinguishability of $d$-neighborhoods, in some sense, in terms of isomorphisms? This would mean that, unlike the game-based framework, we had an alternative to the fact that, for example, FO and $\mathrm{FO}(\mathbf{Q}_p)$ are not Hanf-local under their games.

The new perspective that then immediately emerges from the foregoing is that rather than simply trying to weaken the indistinguishability of neighborhoods when, for example, we do not have sufficient isomorphic neighborhoods to apply locality techniques, we could move to an ''approximation'' of $\mathbf{STRUCT}[\sigma_n]_{(d,0)}^{\widetilde{\mathcal{T}(\sigma_n)}}$ where we could find enough isomorphic neighborhoods. That is, if $\mathbf{STRUCT}[\sigma_n]_{(d,0)}^{\widetilde{\mathcal{T}(\sigma_n)}}$ does not behave well with respect to some property (such as having sufficient isomorphic neighborhoods), why not move on to a $\mathbf{STRUCT}[\sigma_n]_{(d,0)}^{{\widetilde{\mathcal{T}(\sigma_n)}}}$ ''approximation'' where such bad behavior does not occur? This kind of solution comes very naturally when dealing with model structures, i.e. when defining a model structure on $\mathbf{STRUCT}[\sigma_n]_{(d,0)}^{\widetilde{\mathcal{T}(\sigma_n)}}$, its ''approximation'' is simply $\mathbf{STRUCT}[\sigma_n]_{(d,0)}^{{\widetilde{\mathcal{T}(\sigma_n)}}_{\mathrm{core}}}$ (homotopy category of the model structure $\mathbb{M}$), where all weak equivalences (in this case, $k$-logical equivalences, for every primitive-positive sentence of quantifier-rank $k$) become isomorphisms. 

Consider the following diagram: 

\begin{equation}
\xymatrix{& N_d^{\mathfrak{A}}(\vec{a}) \ar[d]_f \ar[r]^{\equiv_k} & N_d^{\mathfrak{B}}(\vec{b}) \ar[d]^g \\ & N_d^{\mathfrak{A}_{\mathfrak{C}}}(\vec{a})  \ar[r]_{\cong} & N_d^{\mathfrak{B}_{\mathfrak{C}}}(\vec{b}),}
\end{equation}

\noindent where $N_d^{\mathfrak{A}}(\vec{a})$ is $\equiv_k$-equivalent to $N_d^{\mathfrak{B}}(\vec{b})$, for primitive-positive sentences of quantifier-rank $k$ (and therefore, by Theorem 1, such $d$-neighborhoods are $k$-homotopically equivalent, i.e. $N_d^{\mathfrak{A}}(\vec{a}) \sim_k N_d^{\mathfrak{B}}(\vec{b})$). 

\noindent $N_d^{\mathfrak{A}_{\mathfrak{C}}}(\vec{a})$ is the core of $N_d^{\mathfrak{A}}(\vec{a})$, $N_d^{\mathfrak{A}_{\mathfrak{C}}}(\vec{a}) \cong N_d^{\mathfrak{B}_{\mathfrak{C}}}(\vec{b})$ is the morphism $N_d^{\mathfrak {A}}(\vec{a}) \sim_k N_d^{\mathfrak{B}}(\vec{b})$ formally inverted in the homotopy category $\mathbf{STRUCT}[\sigma_n]_{{(d,0)}^{\widetilde{\mathcal{T}(\sigma_n)}}_{\mathrm{core}}}$, and $f$, $g$ are homotopic equivalences.

Now, note that if $f$, $g$ are homotopic equivalences, then $f$, $g$ restrict to $k$-homotopic equivalences. So we have the following diagram

\begin{equation}
\xymatrix{& N_d^{\mathfrak{A}}(\vec{a}) \ar[d]_{f_k} \ar[r]^{\equiv_k} & N_d^{\mathfrak{B}}(\vec{b}) \ar[d]^{g_k} \\ & N_d^{\mathfrak{A}_{\mathfrak{C},k}}(\vec{a})  \ar[r]_{\cong} & N_d^{\mathfrak{B}_{\mathfrak{C},k}}(\vec{b}),}
\end{equation}

\noindent where $N_d^{\mathfrak{A}_{\mathfrak{C}, k}}(\vec{a})$ is the $k$-core of $N_d^{\mathfrak{A}}(\vec{a})$; which, by Theorem 1, for every primitive-positive sentence of quantifier-rank $k$, we have

\begin{equation}
\xymatrix{& N_d^{\mathfrak{A}}(\vec{a}) \ar[d]_{\equiv_k} \ar[r]^{\equiv_k} & N_d^{\mathfrak{B}}(\vec{b}) \ar[d]^{\equiv_k} \\ & N_d^{\mathfrak{A}_{\mathfrak{C},k}}(\vec{a})  \ar[r]_{\cong} & N_d^{\mathfrak{B}_{\mathfrak{C},k}}(\vec{b}).}
\end{equation}

What diagram (3) tells us is that for every primitive-positive sentence $\Phi$ of quantifier-rank $k$, $N_d^{\mathfrak{A}}(\vec{a}) \models \Phi \Leftrightarrow N_d^{\mathfrak{B}}(\vec{b}) \models \Phi \Leftrightarrow N_d^{\mathfrak{X}_{\mathfrak{C}, k}}(\vec{x}) \models \Phi$, where $N_d^{\mathfrak{X}_{\mathfrak{C}, k}}(\vec{x})$ is a $k$-core isomorphic to $N_d^{\mathfrak{A}_{\mathfrak{C}, k}}(\vec{a})$ and $N_d^{\mathfrak{B}_{\mathfrak{C}, k}}(\vec{b})$. Therefore, what diagram (3) tells us is that for every primitive-positive sentence $\Phi$ of quantifier-rank $k$, we have 

$$N_d^{\mathfrak{A}}(\vec{a}) \equiv_k N_d^{\mathfrak{B}}(\vec{b}) \Leftrightarrow (N_d^{\mathfrak{A}}(\vec{a}) \equiv_k N_d^{\mathfrak{A}_{\mathfrak{C}, k}}(\vec{a})) \wedge (N_d^{\mathfrak {B}}(\vec{b}) \equiv_k N_d^{\mathfrak{B}_{\mathfrak{C}, k}}(\vec{b})).$$

Thus, the definition of Hanf/Gaifman-localities under $\sim_k$-equivalence (and, therefore, by Theorem 1, under $\equiv_k$-equivalence, for every primitive-positive sentence of quantifier-rank $k$) recovers locality under isomorphism in the homotopy category of $\mathbf{STRUCT}[\sigma_n]_{(d, 0)}^{\widetilde{\mathcal{T}(\sigma_n)}}$:

\begin{dfn}[Gaifman-locality under $(\sim_k,\cong)$-equivalence]
	
	An $m$-ary query $Q$, $m > 0$, on $\sigma_n$-structures, is called \emph{Gaifman-local under $(\sim_k,\cong)$-equivalence} if there exists a number $d \geq 0$ such that for
	every $\sigma$-structure $\mathfrak{A}$ and every $\vec{a}_1, \vec{a}_2 \in A^m$,

		$$N_d^{\mathfrak{A}}(\vec{a}_1) \sim_k N_d^{\mathfrak{A}}(\vec{a}_2) \textrm{ implies } (\vec{a}_1 \in Q(\mathfrak{A}) \Leftrightarrow \vec{a}_2 \in Q(\mathfrak{A})) \Leftrightarrow N_d^{\mathfrak{A}_{\mathfrak{C},k}}(\vec{a}) \cong N_d^{\mathfrak{B}_{\mathfrak{C},k}}(\vec{b})$$ 
		
		$$\textrm{ implies } (\vec{a}_1 \in Q(\mathfrak{A}) \Leftrightarrow \vec{a}_2 \in Q(\mathfrak{A})).$$
		
	\noindent The minimum d for which the above condition holds is called the locality rank of \emph{Gaifman-locality rank of $Q$ under $(\sim_k,\cong)$-equivalence}, and is denoted by $\mathsf{lr}_{(\sim_k,\cong)}(Q)$.
	
\end{dfn}

For Hanf-locality, I define the following:

\begin{dfn}
	
	Consider $\vec{a} \in A^n$ e $\vec{b} \in B^n$. So two $\sigma$-structures $\mathfrak{A},\mathfrak{B}$ are $d$-equivalents under $(\sim_k,\cong)$-equivalence if, and only if, there is a bijection $f: A \rightarrow B$ such that $N_d^{\mathfrak{A}}(\vec{a}c) \sim_k N_d^{\mathfrak{B}}(\vec{b}f(c))$, for every $c \in A$ if, and only if, $N_d^{\mathfrak{A}_{\mathfrak{C},k}}(\vec{a}c) \cong N_d^{\mathfrak{B}_{\mathfrak{C},k}}(\vec{b}f(c))$. When this occurs, I denote this fact by $(\mathfrak{A},\vec{a}) \leftrightarrows_{d,(\sim_k,\cong)} (\mathfrak{B},\vec{b})$. If $n = 0$, we have to satisfy only the condition that there is a bijection $f: A \rightarrow B$ tal que $N_d^{\mathfrak{A}}(c) \sim_k N_d^{\mathfrak{B}}(f(c))$, for every $c \in A$ if, and only if, $N_d^{\mathfrak{A}_{\mathfrak{C},k}}(c) \cong N_d^{\mathfrak{B}_{\mathfrak{C},k}}(f(c))$, and the denotation is simply $\mathfrak{A} \leftrightarrows_{d,(\sim_k,\cong)} \mathfrak{B}$.
	
\end{dfn}

\begin{dfn}[Hanf-locality under $(\sim_k,\cong)$-equivalence]
	
	An $m$-ary query $Q$ on $\sigma$-structures is \emph{Hanf-local under $(\sim_k,\cong)$-equivalence} if there exists a number $d \geq 0$ such that for every $\mathfrak{A}, \mathfrak{B} \in \mathrm{STRUCT}[\sigma]$,
	
	$$(\mathfrak{A},\vec{a}) \leftrightarrows_{d,(\sim_k,\cong)} (\mathfrak{B},\vec{b}) \textrm{ implies } (\vec{a} \in Q(\mathfrak{A}) \Leftrightarrow \vec{b} \in Q(\mathfrak{B})).$$
	
	\noindent The smallest $d$ for which the above condition holds is called the \emph{Gaifman-locality rank of $Q$ under $(\sim_k,\cong)$-equivalence}, and is denoted by $\mathsf{hlr}_{(\sim_k,\cong)}(Q)$.
	
\end{dfn}

I denote this interaction between $\sim_k$-equivalence in the model category $\mathcal{C}$ and isomorphisms in the homotopy category of $\mathcal{C}$ by \emph{$\mathcal{C}_{\leftrightarrows-approximation}$}.

It is also possible to focus on the definition of model structures over $\mathbf{STRUCT}[\sigma_n]_{(d,0)}^{\widetilde{\mathcal{T}(\sigma_n)}}$ in order to investigate the properties of their homotopic equivalences with respect to locality. For example, there are three trivial model structures over $\mathbf{STRUCT}[\sigma_n]_{(d,0)}^{\widetilde{\mathcal{T}(\sigma_n)}}$, where the choice of subcategories of fibrations, cofibrations, and weak equivalences are reduced to $\mathbf{STRUCT}[\sigma_n]_{(d,0)}^{\widetilde{\mathcal{T}(\sigma_n)}}$ and $\mathbf{STRUCT}[\sigma_n]^{\widetilde{\mathcal{T}(\sigma_n)}}_{{(d,0)}_{\mathrm{iso}}}$ (its restriction to isomorphisms). In the case where we have a model structure $\mathbb{M}$ over $\mathbf{STRUCT}[\sigma_n]_{(d,0)}^{\widetilde{\mathcal{T}(\sigma_n)}}$ whose subcategory of weak equivalences is $\mathbf{STRUCT}[\sigma_n]^{\widetilde{\mathcal{T}(\sigma_n)}}_{{(d,0)}_{\mathrm{iso}}}$, trivially follows that the locality under weak equivalences is only the usual locality under isomorphisms. Thus, it is possible to classify and investigate locality under different equivalences by investigating possible model structures over $\mathbf{STRUCT}[\sigma_n]_{(d, 0)}^{\widetilde{\mathcal{T}(\sigma_n)}}$.

Given the above idea, I propose a general definition of locality for  Quillen model category-based frameworks. To do this, let us now return to diagrams (1), (2) and (3). Given that $ N_d^{\mathfrak{A}_{\mathfrak{C}, k}}(\vec{a}) \cong N_d^{\mathfrak{B}_{\mathfrak{C}, k}}(\vec{b})$ is an isomorphism, we can consider a single object (up to isomorphism), say $N_d^{\mathfrak{X}_{\mathfrak{C}, k}}(\vec{x})$. What gives us

\begin{equation}
\xymatrix{& N_d^{\mathfrak{A}}(\vec{a}) \ar[dr]_{\sim_k} \ar[r]^{\sim_k} & N_d^{\mathfrak{B}}(\vec{b}) \ar[d]^{\sim_k} \\ &  & N_d^{\mathfrak{X}_{\mathfrak{C},k}}(\vec{x}).}
\end{equation}

First, note that in a given category $\mathcal{C}$, by the Yoneda Lemma, a morphism $f: A \rightarrow B$ is an isomorphism precisely in the case that for every object $X$,

$$\mathrm{Hom}_{\mathcal{C}}(f,X): \mathrm{Hom}_{\mathcal{C}}(B,X) \rightarrow \mathrm{Hom}_{\mathcal{C}}(A,X)$$

\noindent it is a bijection. We can then use this for the case where we have a category $\mathcal{C}$ and a subcategory $\mathcal{E}$ such that morphisms in $\mathcal{E}$ are some kind of equivalence we want treat as isomorphisms. For the context we are dealing with, $\mathcal{E}$ has homotopic equivalences as morphisms. Thus, it is possible to define certain $\mathcal{C}$-objects called $\mathcal{E}$-locals as follows: An object $X$ of $\mathcal{C}$  is said to be $\mathcal{E}$-local if any morphism $f: A \rightarrow B$ in $\mathcal{E}$ induces a bijection

$$f^*: \mathrm{Hom}_{\mathcal{C}}(B,X) \cong \mathrm{Hom}_{\mathcal{C}}(A,X).$$

Note that diagram (4) shows us that 
$$N_d^{\mathfrak{A}}(\vec{a}) \sim_k N_d^{\mathfrak{B}}(\vec{b}) \Leftrightarrow (N_d^{\mathfrak{A}}(\vec{a}) \sim_k N_d^{\mathfrak{X}_{\mathfrak{C},k}}(\vec{x})) \wedge (N_d^{\mathfrak{B}}(\vec{b}) \sim_k N_d^{\mathfrak{X}_{\mathfrak{C},k}}(\vec{x})).$$ 
But, this is exactly the same as saying that there is a $d$-neighborhood $N_d^{\mathfrak{X}_{\mathfrak{C},k}}(\vec{x})$ in $\mathbf{STRUCT}[\sigma_n]_{(d,0)}^{\widetilde{\mathcal{T}(\sigma_n)}}$ such that a $k$-homotopic equivalence $f_k : N_d^{\mathfrak{A}}(\vec{a}) \sim_k N_d^{\mathfrak{B}}(\vec{b})$ in a subcategory of $\mathbf{STRUCT}[\sigma_n]_{(d,0)}^{\widetilde{\mathcal{T}(\sigma_n)}}$ induces a bijection

$$f^*:\mathrm{Hom}_{\mathbf{STRUCT}[\sigma_n]_{(d,0)}^{\widetilde{\mathcal{T}(\sigma_n)}}}(N_d^{\mathfrak{B}}(\vec{b}),N_d^{\mathfrak{X}_{\mathfrak{C},k}}(\vec{x})) \cong \mathrm{Hom}_{\mathbf{STRUCT}[\sigma_n]_{(d,0)}^{\widetilde{\mathcal{T}(\sigma_n)}}}(N_d^{\mathfrak{A}}(\vec{a}),N_d^{\mathfrak{X}_{\mathfrak{C},k}}(\vec{x})).$$

So for a model structure $\mathbb{M}$ over $\mathbf{STRUCT}[\sigma_n]_{(d,0)}^{\widetilde{\mathcal{T}(\sigma_n)}}$, and the subcategory $$\mathbf{STRUCT}[\sigma_n]^{\widetilde{\mathcal{T}(\sigma_n)}}_{{(d,0)}_{{we}}} \subset \mathbf{STRUCT}[\sigma_n]_{(d,0)}^{\widetilde{\mathcal{T}(\sigma_n)}}$$ of weak equivalences $\sim$, I propose the following general definition:

\begin{dfn}[Gaifman-locality under Hom-isomorphisms of $\sim$-equivalence]
	
	An m-ary query $Q$, $m > 0$, on $\sigma$-structures, is called \emph{Gaifman-local under Hom-isomorphisms of $\sim$-equivalence} if there exists a number $d \geq 0$ such that for every $\sigma$-structure $\mathfrak{A}$, and every $\vec{a}_1,\vec{a}_2 \in A^m$, there is an object $\mathbf{STRUCT}[\sigma_n]^{\widetilde{\mathcal{T}(\sigma_n)}}_{{(d,0)}_{{we}}}$-local $N_d^{\mathfrak{X}_{\mathfrak{C}}}(\vec{x})$ of $\mathbf{STRUCT}[\sigma_n]_{(d,0)}^{\widetilde{\mathcal{T}(\sigma_n)}}$, and an equivalence $f:N_d^{\mathfrak{A}}(\vec{a}) \sim N_d^{\mathfrak{B}}(\vec{b})$ in $\mathbf{STRUCT}[\sigma_n]^{\widetilde{\mathcal{T}(\sigma_n)}}_{{(d,0)}_{{we}}}$, inducing a bijection
	
	$$f^*: \mathrm{Hom}_{\mathbf{STRUCT}[\sigma_n]_{(d,0)}^{\widetilde{\mathcal{T}(\sigma_n)}}}(N_d^{\mathfrak{B}}(\vec{b}),N_d^{\mathfrak{X}_{\mathfrak{C}}}(\vec{x})) \cong \mathrm{Hom}_{\mathbf{STRUCT}[\sigma_n]_{(d,0)}^{\widetilde{\mathcal{T}(\sigma_n)}}}(N_d^{\mathfrak{A}}(\vec{a}),N_d^{\mathfrak{X}_{\mathfrak{C}}}(\vec{x})),$$   
	
	\noindent that implies
	$$(\vec{a}_1 \in Q(\mathfrak{A}) \Leftrightarrow \vec{a}_2 \in Q(\mathfrak{A})).$$
	\noindent The minimum $d$ for which the above condition holds is called the \emph{Gaifman-locality rank of $Q$ under Hom-isomorphisms of $\sim$-equivalence}, and is denoted by $\mathsf{lr}_{\mathrm{iso}(\mathbb{M}_{we})}(Q)$.
	
\end{dfn}

\begin{dfn}[Hanf-locality under Hom-isomorphisms of $\sim$-equivalence]

An $m$-ary query $Q$ on $\sigma$-structures is \emph{Hanf-local under Hom-isomorphisms of $\sim$-equivalence} if there exists a number $d \geq 0$ such that for every $\mathfrak{A}, \mathfrak{B} \in \mathbf{STRUCT}[\sigma]$, there is a bijection $f: A \rightarrow B$ such that $h : N_d^{\mathfrak{A}}(\vec{a}c) \sim N_d^{\mathfrak{B}}(\vec{b}f(c))$ in $\mathbf{STRUCT}[\sigma_n]^{\widetilde{\mathcal{T}(\sigma_n)}}_{{(d,0)}_{{we}}}$, and an object $\mathbf{STRUCT}[\sigma_n]^{\widetilde{\mathcal{T}(\sigma_n)}}_{{(d,0)}_{{we}}}$-local $X$ of $\mathbf{STRUCT}[\sigma_n]_{(d,0)}^{\widetilde{\mathcal{T}(\sigma_n)}}$, inducing a bijection
    
 $$f^*: \mathrm{Hom}_{\mathbf{STRUCT}[\sigma_n]_{(d,0)}^{\widetilde{\mathcal{T}(\sigma_n)}}}(N_d^{\mathfrak{B}}(\vec{b}f(c)),X) \cong \mathrm{Hom}_{\mathbf{STRUCT}[\sigma_n]_{(d,0)}^{\widetilde{\mathcal{T}(\sigma_n)}}}(N_d^{\mathfrak{A}}(\vec{a}c),X),$$   
 
 \noindent that implies
 
 $$(\vec{a} \in Q(\mathfrak{A}) \Leftrightarrow \vec{b} \in Q(\mathfrak{B})).$$
 
 The smallest $d$ for which the above condition holds is called the \emph{rank da Hanf-localidade de $Q$ under Hom-isomorphisms of $\sim$-equivalence}, and is denoted by $\mathsf{hlr}_{\mathrm{iso}(\mathbb{M}_{we})}(Q)$.

\end{dfn}

Definitions 8 and 9 are interesting not only because they are natural isomorphisms, but also because they allow to describe the concept of localization within the scope of locality of logics, as it usually comes up in topology. 


\section{Modules and Pre-Triangulated Category}

In this section we briefly present how the homotopy category of any model category is naturally a $\mathrm{Ho(SSet)}$-module category. The main point of this section is the fact that the homotopy category of a pointed model category is a pre-triangulated category.

\subsection{Modules}

\begin{dfn}
Let $\mathcal{C}$, $\mathcal{D}$ and $\mathcal{E}$ be categories. An \emph{adjunction of two variables} consists of functors
$$- \otimes -: \mathcal{C} \times \mathcal{D} \rightarrow \mathcal{E}$$
$$(-)^{(-)}: \mathcal{D}^{\mathrm{op}} \times \mathcal{E} \rightarrow \mathcal{C}$$
$$\mathrm{map}(-,-): \mathcal{C}^{\mathrm{op}} \times \mathcal{E} \rightarrow \mathcal{D}$$
\noindent satisfying the usual adjointness conditions. See \emph{[12, Definition 4.1.12]}.
\end{dfn}{}

If $\mathcal{C}$, $\mathcal{D}$ and $\mathcal{E}$ are model categories, we would like to know how an adjunction of two variables can be compatible with the respective model structures. This occurs as follows.

\begin{dfn}
Now let $\mathcal{C}$, $\mathcal{D}$ and $\mathcal{E}$ be model categories. A \emph{Quillen adjunction of two variables} is an adjunction of two variables such that: If $f: U \rightarrow V$ is a cofibration in $\mathcal{C}$ and $g: W \rightarrow X$ is a cofibration in $\mathcal{D}$, then the induced pushout-product map
$$f \msquare g: (U \otimes X) \coprod_{U \otimes W} (V \otimes W) \rightarrow V \otimes X$$
\noindent is a cofibration in $\mathcal{E}$. Furthermore, the map $f \msquare g$ must be a trivial cofibration if either of $f$ or $g$ is.
\end{dfn}{}

\begin{dfn}
Let $\mathcal{D}$ be a closed symmetric monoidal category with product $\times$ and unit $S$. A category $\mathcal{M}$ is a closed $\mathcal{D}$-module category if it has an adjunction of two variables
$$(- \otimes -,(-)^{(-)},\mathrm{map}(-,-)): \mathcal{M} \times \mathcal{D} \rightarrow \mathcal{M}$$
\noindent together with natural associativity isomorphisms
$$(X \otimes D) \otimes E \rightarrow X \otimes (D \times E)$$
\noindent and natural unit isomorphisms
$$X \otimes S \rightarrow X.$$
\noindent These isomorphisms have to satisfy some standard coherence conditions.
That is, the pentagonal diagram describing fourfold associativity must commute,
as must the triangle relating the two ways to obtain $X \otimes D$ from $X \otimes (S \times D)$.
\end{dfn}{}

\begin{dfn}
Let $\mathcal{D}$ be a closed symmetric monoidal model category. A
model category $\mathcal{M}$ is a $\mathcal{D}$-model category if it is a $\mathcal{D}$-module category in the sense of Definition 11 satisfying the following.
\begin{itemize}
    \item $- \otimes -$ is a Quillen bifunctor.
    \item Let $QS \rightarrow S$ be the cofibrant replacement of the unit in $\mathcal{D}$ and let $X \in \mathcal{M}$ be cofibrant. Then
    $$X \otimes QS \rightarrow X \otimes S$$
    \noindent is a weak equivalence in $\mathcal{M}$.
\end{itemize}{}
\end{dfn}{}

\subsection{Framings}

We will now talk about the notions of cosimplicial and simplicial frames. Such notions appear to solve the following problem: suppose one is studying a model category $\mathcal{C}$ that is not necessarily simplicial, one would still like to have a reasonable substitute for tensoring with simplicial sets or for mapping spaces. Framings provide such a generalisation.

The idea is to take an object $A \in \mathcal{C}$ and then apply a particular
cofibrant (respectively fibrant) replacement. Here, $A \in \mathcal{C}$ is being considered as a constant cosimplicial (or simplicial object) in $\mathcal{C}$.

The resulting cosimplicial or simplicial objects can then be used to solve the problem. However, this will not make $\mathcal{C}$ a simplicial model category. But it can at least ensure that the homotopy category $\mathrm{Ho}(\mathcal{C})$ is a closed Ho(SSet)-module (Where SSet denotes the category of simplicial sets). For more details on framings see [Hov99, Chapter 5].

Let $\mathcal{C}$ be a category. By $\mathcal{C}^{\Delta}$ we denote the category of cosimplicial objects in $\mathcal{C}$. The standard model structure for this category is the Reedy model structure, which is described in [Hov99, Section 5.1].We started with a very useful proposition.

\begin{pps}
Suppose $\mathcal{C}$ is a category with all small colimits. Then the
category $\mathcal{C}^{\Delta}$ is equivalent to the category of adjunctions $\mathrm{SSet} \rightleftarrows \mathcal{C}$. We denote the image of $A^{\bullet} \in \mathcal{C}^{\Delta}$ under this equivalence by $(A^{\bullet} \otimes -, \mathcal{C}(A^{\bullet},-),\varphi): \mathrm{SSet} \rightarrow \mathcal{C}$.
\end{pps}{}

\begin{proof}
(12, Proposition 3.1.5)
\end{proof}{}

\begin{rem} Dually, if $\mathcal{C}$ has all small limits, there is an equivalence of categories between $\mathcal{C}^{\Delta^{\mathrm{op}}}$ (the category of simplicial objects) and adjunctions $\mathrm{SSet}^{\mathrm{op}} \leftrightarrows \mathcal{C}$. We denote the image of a simplicial object $A_{\bullet}$ by $(\mathrm{Hom}(-,A_{\bullet}), \mathcal{C}(-,A_{\bullet}),\varphi)$. We might also write $\mathrm{Hom}(-,A_{\bullet}) = A_{\bullet}^{(-)}$.
\end{rem}

In addition, we have:

\begin{itemize}
    \item $A^{\bullet} \otimes \Delta[n] = A^{\bullet}[n]$,
    \item $A^{\bullet} \otimes \partial \Delta[n] \rightarrow A^{\bullet} \otimes \Delta[n]$ is the $n^{th}$ latching map of $A^{\bullet}$ [Hir03, Proposition 16.3.8].
    \item $A^{\bullet} \otimes -$ preserves colimits.
\end{itemize}{}

Dually, we have:

\begin{itemize}
    \item $A_{\bullet}^{\Delta[n]} = A_{\bullet}[n]$.
    \item $A_{\bullet}^{\Delta[n]} \rightarrow A_{\bullet}^{\partial \Delta[n]}$ is the $n^{th}$ matching map of $A_{\bullet}$ [Hir03, Proposition 16.3.8].
    \item $A_{\bullet}^{(-)}$ takes limits of SSet to colimits of $\mathcal{C}$.
\end{itemize}{}

\begin{dfn}
If $\mathcal{C}$ is a model category, we say that an object $A^{\bullet} \in \mathcal{C}^{\Delta}$ is a \emph{cosimplicial frame} if
$$A^{\bullet} \otimes - :\mathrm{SSet} \rightleftarrows \mathcal{C}: \mathrm{map}_l(A^{\bullet},-)$$
\noindent is a Quillen adjunction.

An object $A_{\bullet} \in \mathcal{C}^{\Delta^{\mathrm{op}}}$ is a \emph{simplicial frame} if
$$A_{\bullet}^{(-)} :\mathrm{SSet}^{\mathrm{op}} \leftrightarrows \mathcal{C}: \mathrm{map}_r(-,A_{\bullet})$$
\noindent is a Quillen adjunction.
\end{dfn}{}

Cosimplicial frames can be characterised as follows.

\begin{pps}
A cosimplicial object $A^{\bullet} \in \mathcal{C}^{\Delta}$ is a cosimplicial frame if and only if $A^{\bullet}$ is cofibrant and the structure maps $A^{\bullet}[n] \rightarrow A^{\bullet}[0]$ are weak equivalences for $n \geq 0$.

\end{pps}{}

\begin{proof}
The ingredients to the proof can be found in [12, Proposition 3.6.8, Example 5.2.4, Theorem 5.2.5, Proposition 5.4.1]. 
\end{proof}{}

The case for simplicial frames is dual.

\begin{thm}[Hovey]
There exists a functor $\mathcal{C} \rightarrow \mathcal{C}^{\Delta}$ such that the image $A^{*}$ of any cofibrant $A \in \mathcal{C}$ under this functor is a cosimplicial frame with $A^{*}[0] \cong A$.

There exists a functor $\mathcal{C} \rightarrow \mathcal{C}^{\Delta^{\mathrm{op}}}$ such that the image $A_{*}$ of any fibrant $A \in \mathcal{C}$ under this functor is a simplicial frame with $A_{*}[0] \cong A$.
\end{thm}{}

\begin{dfn}
A functor $A \mapsto A^{*}$ together with a functor $A \mapsto A_{*}$ satisfying
the conditions of Theorem 1 is called a \emph{framing} of $\mathcal{C}$.
\end{dfn}{}

Together with the framing functors $A \mapsto A^{*}$ and $A \mapsto A_{*}$ of Theorem 3 one obtains bifunctors

\begin{itemize}
    \item $- \otimes -: \mathcal{C} \times \mathrm{SSet} \rightarrow \mathcal{C}, \textrm{ } (A,K) \mapsto A^{*} \otimes K$.
    \item $\mathrm{map}_l(-,-): \mathcal{C}^{\mathrm{op}} \times \mathcal{C} \rightarrow \mathrm{SSet}, \textrm{ } (A,B) \mapsto \mathrm{map}_l(A^{*},B)$.
    \item $(-)^{(-)}:\mathrm{SSet^{op}} \times \mathcal{C} \rightarrow \mathcal{C}, \textrm{ } (A,K) \mapsto A_{*}^K$.
    \item $\mathrm{map}_r(-,-): \mathcal{C}^{\mathrm{op}} \times \mathcal{C} \rightarrow \mathrm{SSet}, \textrm{ } (A,B) \mapsto \mathrm{map}_r(A,B_{*})$.
\end{itemize}{}

This does not make $\mathcal{C}$ into a simplicial model category. But, Hovey shows in [12, Theorem 5.4.9] that

$$- \otimes -: \mathcal{C} \times \mathrm{SSet} \rightarrow \mathcal{C}$$
\noindent and
$$(-)^{(-)}:\mathrm{SSet} \times \mathcal{C}^{\mathrm{op}} \rightarrow \mathcal{C}^{\mathrm{op}}$$

\noindent (with the opposite model structure) have total left derived functors.

By [12, Proposition 5.4.7], we know that the two right adjoints $\mathrm{map}_l$ and $\mathrm{map}_r$ only agree up to a zig-zag of weak equivalences in $\mathcal{C}$. However, this means the right derived mapping spaces $R \mathrm{map}_l$ and $R \mathrm{map}_r$ agree. However, this gives us an adjunction of two variables 

$$(- \otimes^L -, R(-)^{(-)}, R\mathrm{map}(-,-)):\mathrm{Ho}(\mathcal{C}) \times \mathrm{Ho}(\mathrm{SSet}) \rightarrow \mathrm{Ho}(\mathcal{C}).$$

Note that the functor $- \otimes -$ is not, in general, associative. But, when we move to the homotopy category, this problem is solved. Hovey details the construction of a particular associativity weak equivalence, and, thus, we have [12, Theorem 5.5.3]:

\begin{thm}[Hovey]
The framing functor of Theorem 3 makes $\mathrm{Ho}(\mathcal{C})$ into a closed $\mathrm{Ho}(\mathrm{SSet})$-module category.
\end{thm}{}

In particular, this result is also valid in the following configuration: if $\mathcal{C}$ is a pointed model category, then $\mathrm{Ho}(\mathcal{C})$ is a
closed $\mathrm{Ho}(\mathrm{SSet}_{*})$-module, where $\mathrm{SSet}_{*}$ denotes the category of pointed simplicial sets.

\subsection{Pre-Triangulated Categories}

As seen above, the homotopy category of a model category is naturally
a closed $\mathrm{Ho(SSet)}$-module, and the homotopy category of a pointed model category is naturally a closed $\mathrm{Ho(SSet_*)}$-module. With that information, we can work with suspension and loop functors. These exist in any closed $\mathrm{Ho(SSet_*)}$-module, but there are a number of results specific to the homotopy category of a pointed model category. We can then use the closed action of $\mathrm{Ho(SSet_*)}$ on $\mathrm{Ho}(\mathcal{C})$ (given in [12, Section 5.7]) to define suspension and loop functors as follows.

\begin{dfn}
Suppose $\mathcal{C}$ is a pointed model category. The \emph{suspension
functor} $\Sigma: \mathrm{Ho}(\mathcal{C}) \rightarrow \mathrm{Ho}(\mathcal{C})$ is the functor $X \mapsto X \wedge^L S^1$ defined by the closed action
of $\mathrm{Ho(SSet_*)}$ on $\mathrm{Ho}(\mathcal{C})$. Dually, the \emph{loop functor} $\Omega: \mathrm{Ho}(\mathcal{C}) \rightarrow \mathrm{Ho}(\mathcal{C})$ is the functor $X \mapsto R\mathrm{Hom}_*(S^1,X)$.
\end{dfn}{}

The suspension functor is of course left adjoint to the loop functor. 

With that, and a few more definitions (see [12, Section 6.1 and Section 6.2]), we can show that there is a natural coaction in $\mathrm{Ho}(\mathcal{C})$ of the cogroup $\Sigma A$ on the cofiber of a cofibration of cofibrant objects $A \rightarrow B$ in a pointed model category $\mathcal{C}$. This allows to define cofiber sequences, and, by duality, fiber sequences. With the construction of a coaction, Theorem 6.2.1 [12] and some other things, we can define the following.

\begin{dfn}
Suppose $\mathcal{C}$ is a pointed model category. A \emph{cofiber sequence}
in $\mathrm{Ho}\mathcal{C}$ is a diagram $X \rightarrow Y \rightarrow Z$ in $\mathrm{Ho}\mathcal{C}$ together with a right coaction of $\Sigma X$
on Z which is isomorphic in $\mathrm{Ho}\mathcal{C}$ to a diagram of the form $A \xrightarrow{f} B \xrightarrow{g} C$ where $f$ is a cofibration of cofibrant objects in $\mathcal{C}$ with cofiber $g$ and where $\mathcal{C}$ has the right $\Sigma A$-coaction given by Theorem 6.2.1. Dually, a \emph{fiber sequence} is a diagram $X \rightarrow Y \rightarrow Z$ together with a right action of $\Omega Z$ on $X$ which is isomorphic to a diagram $F \xrightarrow{i} E \xrightarrow{p} B$ where $p$ is a fibration of fibrant objects with fiber $i$ and where $F$ has the right $\Omega B$-action given by Theorem 6.2.1.
\end{dfn}{}

A cofiber sequence has associated to it a boundary map, the definition is as follows.

\begin{dfn}
Suppose $\mathcal{C}$ is a pointed model category, and $X \xrightarrow{f} Y \xrightarrow{g} Z$ is a cofiber sequence in $\mathrm{Ho}\mathcal{C}$. The \emph{boundary map} is the map $\partial: Z \rightarrow \Sigma X$ in $\mathrm{Ho}\mathcal{C}$ which is the composite
$$Z \rightarrow Z \coprod \Sigma X \xrightarrow{0 \times 1}  \Sigma X$$
\noindent where the first map is the coaction.
Dually, if $X \xrightarrow{f} Y \xrightarrow{g} Z$ is a
fiber sequence, the \emph{boundary map} is the map $\partial: \Omega Z \rightarrow X$ which is the composite
$$\Omega Z \xrightarrow{0 \times 1} X \times \Omega Z \rightarrow Z.$$
\end{dfn}{}

For more details, see [12, Section 6.2].

The point is that cofiber and fiber sequences have certain properties (see [12, Section 6.3 and Section 6.4] that can be abstracted to define the notion of pre-triangulation, which in turn leads us to the notion of pre-triangulated category. To be more exact, suppose $\mathcal{S}$ is a nontrivial (right) closed $\mathrm{Ho(SSet_)}$-module. A pre-triangulation on $\mathcal{S}$ is a collection of cofiber sequences and fiber sequences satisfying certain conditions (see [12, Definition 6.5.1]). A pre-triangulated category is then a nontrivial closed $\mathrm{Ho(SSet_*)}$-module $\mathcal{S}$ with all small coproducts and products, together with a pre-triangulation
on $\mathcal{S}$. This shows that the homotopy category of a pointed
model category is a pre-triangulated category.

\section{Pre-Triangulated Categories and the Category of Structures}

From Theorem 2, we immediately know that the homotopy category of $\mathbf{STRUCT}[\sigma]$ is a $\mathrm{Ho(SSet)}$-module category. It is easy to see that $\mathbf{STRUCT}[\sigma]$ is not a pointed model category. However, as we will see in the next subsection, using over and under categories it is possible to build pointed categories from categories that are not pointed categories. Furthermore, if $\mathcal{C}$ is a model category, the pointed category built from $\mathcal{C}$ inherits the model structure.

\subsection{General Pointed Category}

\begin{dfn}
If $A$ be an object in $\mathcal{C}$ , the categories under $A$ and over $A$ will be denoted by $\mathcal{C}^A$ , $\mathcal{C}_A$, respectively. Objects and morphisms in $\mathcal{C}_A$ will be underlined, $\underline{f}: \underline{X} \rightarrow \underline{Y}$ , and for the category over $A$ the overlined notation $\overline{f}: \overline{X} \rightarrow \overline{Y}$ will be used. A category $\mathcal{C}$ is said to be pointed if there exist initial and final objects and they are isomorphic. This object is usually denoted by $*$ and
it is called the zero object.
\end{dfn}{}

In addition, for every category $\mathcal{C}$ one has the following properties:

\begin{enumerate}
    \item $\mathcal{C}^A$ always has initial object $\overline{A} = \mathrm{id}_A: A \rightarrow A$.
    \item If $A$ is the initial object of $\mathcal{C}$ , then $\mathcal{C}^A \cong \mathcal{C}$.
    \item If $A$ is the terminal object of $\mathcal{C}$ , then $\mathcal{C}^A$ is a pointed category, where $\overline{A}: A \rightarrow A$ is the zero object.
\end{enumerate}{}

And we also have the corresponding dual properties.

\begin{dfn}
If $A$ is any given object in a category $\mathcal{C}$,
$$(\mathcal{C}^A)_{\overline{A}} = (\mathcal{C}_A)^{\overline{A}}$$
\noindent is a pointed category that it will also be denoted by $\mathcal{C}_A^A$. An object in this category $\underline{\overline{X}}: A \xrightarrow{i_X} X \xrightarrow{r_X} A$ is determined by morphisms $i_X$ and $r_X$ in $\mathcal{C}$ such that $r_X i_X = \mathrm{id}_A$. A morphism $\underline{\overline{f}}: \underline{\overline{X}} \rightarrow \underline{\overline{Y}}$ is given by a morphism $f:X \rightarrow Y$ in $\mathcal{C}$ such that $f i_X = i_Y$ and $r_Y f = r_X$. The zero object of $\mathcal{C}_A^A$ is $\underline{\overline{A}}: A \xrightarrow{\mathrm{id}_A} A \xrightarrow{\mathrm{id}_A} A$.  
\end{dfn}{}

For a given object $A$ in a category $\mathcal{C}$ , we can consider the forgetful functors $U: \mathcal{C}^A \rightarrow \mathcal{C}$ and $V: \mathcal{C}_A \rightarrow \mathcal{C}$, given by $U(\overline{X}) = X$, $U(\overline{f}) = f$ and $V(\underline{X}) = X$, $V(\underline{f}) = f$.

If $\mathcal{C}$ has finite coproducts, then $U$ has a left adjoint $A \overline{\sqcup} (-): \mathcal{C} \rightarrow \mathcal{C}^A$ given as
follows. If $X$ is an object in $\mathcal{C}$ , $A \overline{\sqcup} X$ is the canonical morphism $\mathcal{C}$, $A \overline{\sqcup} X: A \rightarrow A \sqcup X$.
Dually, if $\mathcal{C}$ has finite products, then $V$ has a right adjoint $A \underline{\times} (-): \mathcal{C} \rightarrow \mathcal{C}_A$, which assigns to an object $X$ in $\mathcal{C}$ , the object $A \underline{\times} X$ which is the projection $A \times X \rightarrow A$.

\begin{dfn}
Let $\mathcal{C}$ be a model category and let $A$ be an object in $\mathcal{C}$.
A morphism $\overline{f}$ in $\mathcal{C}^A$ is said to be a cofibration, a fibration or a weak equivalence if and only if $U(\overline{f})$ is a cofibration, a fibration or a weak equivalence in $\mathcal{C}$. 

In a dual way, a morphism $\underline{f}$ in $\mathcal{C}_A$ is said to be a cofibration, a fibration or a weak equivalence if and only if $V(\underline{f})$ is a cofibration, a fibration or a weak equivalence in $\mathcal{C}$.
\end{dfn}{}

We refer the reader to [18] for a proof of the following:

\begin{pps}
Let $\mathcal{C}$ be a model category and let $A$ be an object in
$\mathcal{C}$. Then, the categories $\mathcal{C}^A$ and $\mathcal{C}_A$, with the classes of morphisms given above, has the structure of a model category.
\end{pps}{}

Thus, we can apply Proposition 3 to $\mathbf{STRUCT}[\sigma]$, and obtain a pointed model category from $\mathbf{STRUCT}[\sigma]$.  With this, we can use the technical apparatus of pre-triangulated categories within the scope of finite structures. Pre-triangulated categories are the unstable analog of triangulated categories, and it has been investigated and studied by many. See, for example, [16].

The main interest in the investigation of pre-triangulated categories within the scope of finite structures is the relationship between pre-triangulated categories and triangulated categories in this context. To be more exact, a triangulated category is a pre-triangulated category in which the suspension functor $\Sigma$ is an equivalence of categories. A pointed model category is stable if its homotopy category is triangulated. 

If we are able to show that the pointed model category obtained from $\mathbf{STRUCT}[\sigma]$ is stable, we will have a wide range of unexplored territory within the scope of finite structures. The main, and most interesting, is the cohomological territory: triangulated categories admit a notion of cohomology, and every triangulated category has a large supply of cohomological functors. 

\section{Final Considerations}

Throughout this paper I have presented the implications of a Quillen model structures-based framework for locality under logical equivalence. However, one point of my proposal remains problematic. As noted in Theorem 1, $k$-homotopic equivalence of $d$-neighborhoods only implies $k$-logical equivalence for primitive-positive sentences of quantifier-rank $k$. That is, $k$-homotopic equivalence of $d$-neighborhoods does not imply $k$-logical equivalence of $d$-neighborhoods for every sentence of quantifier-rank $k$ (and similarly to Theorem 2).

So my goal in future developments is to extend $k$-homotopic equivalence to imply not only $k$-logical equivalence for primitive-positive sentences of quantifier-rank $k$, but to imply $k$-logical equivalence for every sentence of quantifier-rank $k$. In addition, it is of obvious interest to investigate the behavior of the bi-implication ''$k$-homotopic equivalence $\Leftrightarrow$ $k$-logical equivalence'' in logics other than FO. 

However, as seen, Corollary 1 ensures the bi-implication ''$k$-homotopic equivalence $\Leftrightarrow$ $k$-logical equivalence'' for $k$-extendable structures, which gives us, at least, a partial view of how the implications seen above can work without restriction.

\bibliographystyle{unsrt}  


\end{document}